\documentclass[a4paper,10pt]{amsart}
\usepackage{amssymb} \usepackage{amsfonts}
\usepackage{amsmath,amsthm}
\usepackage{tikz,relsize}
\usepackage[cmtip,arrow]{xy}
\usepackage{pb-diagram,pb-xy}
\newtheoremstyle{mytheoremstyle}{2ex}{2ex}{\it}{0pt}{\bf}{.}{ }{}
\newtheoremstyle{myalgorithmstyle}{2ex}{2ex}{}{0pt}{\rm}{.}{ }{}
\theoremstyle{mytheoremstyle}
\newtheorem{lm}{Lemma}[section]
\newtheorem{theo}{Theorem}[section]
\newtheorem{cor}{Corollary}[section]
\newtheorem{prop}{Proposition}[section]
\theoremstyle{myalgorithmstyle}
\newtheorem{defi}{Definition}[section]
\newtheorem{ex}{Example}[section]

\begin{document}
	
\title[Spin-structures on real Bott manifolds with K\"{a}hler structure]{Spin-structures on real Bott manifolds with K\"{a}hler structure}

\author{A. G\c{a}sior}

\begin{abstract}  Let $M$ be a real Bott manifold with K\"{a}hler structure. Using Ishida characterization \cite{I} we give necessary and sufficient condition for the existence of the Spin-structure on $M$. In proof we use the technic developed in \cite{PS} and characteristic classes.
\end{abstract}
\keywords{real Bott manifolds, Spin structure, K\"{a}hler structure\newline Author is supported by the Polish National Science Center grant DEC-2013/09/B/ST1/04125}

\maketitle

\hskip5mm
 \section{Introduction }

Let $\Gamma$ be a fundamental group of a real Bott manifold $M$. From \cite{KM} we know that $\Gamma$ defines a short exact sequence
\begin{equation}
0\to \mathbb Z^n\xrightarrow{\iota}\Gamma\xrightarrow{\pi} \mathbb Z_2^k\to1.\label{wzor21}
\end{equation}
However above we have a induced holonomy representation $\varrho:\mathbb Z_2^k\to GL(n,\mathbb Z)$
$$\varrho(g)(z)=\iota(\gamma\iota(z)\gamma^{-1})$$
for all $g\in \mathbb Z_2^k$, $\pi(\gamma)=g$, $\gamma\in\Gamma$, $z\in\mathbb Z^n$, where $\varrho(\mathbb Z_2)^k\subset D=\{A=(a_{ij})\in GL(n,\mathbb Z)|a_{ij}=0,i\ne j, a_{ii}=\pm1,1\leq i,j\leq n\}.$

It is well known, (\cite{KM}), that $M$ is determined by a certain matrix $A=[a_{ij}]$, where $a_{ij}\in\mathbb F_2$. We call $A$ a Bott matrix and we shall denote the manifold $M$  by $M(A)$. In \cite{I} Ishida gave the following necessary and sufficient condition for existence of the K\"{a}hler structure on $M(A)$.
\begin{theo}{\rm (\cite{I} Theorem 3.1)}\label{thI}
	Let $A$ be $2n-$dimensional matrix of real Bott manifold $M(A)$. Then the following conditions are equivalent:
	\begin{enumerate}
		\item there exist $n$ subsets $\{j_1,j_{n+1}\}$, $\ldots$, $\{j_n,j_{2n}\}$ of $\{1,2,\ldots,2n\}$ such that
		\begin{enumerate}
			\item	$\coprod_{k=1}^n\{j_k,j_{k+n}\}=\{1,2,\ldots,2n\}$,
			\item $A^{j_k}=A^{j_{k+n}}$ for all $1\leq i<j_k$, where $A^k$ is the $k-$th column of the matrix $A$.
		\end{enumerate}
		\item there exist a K\"{a}hler structure on $M(A)$.	
	\end{enumerate}
\end{theo}
In this note we are going to examine the existence of the Spin structure on a real Bott manifold with a K\"{a}hler structure. We would like to mention that the general condition for existence of Spin structure is considered in (\cite{GS1}, \cite{G1}). However we use different new methods which we find interesting. We use different definition of the real Bott manifold which we introduce  at the section 2. In the section 3 we prove our main result.

\section{New definition of real Bott manifold}

\indent In this section we recall methods introduced in \cite{PS} and developed in \cite{LPPS}.
Let $S^1$ be a unit circle in $\mathbb C$ and we consider authomorphisms $g_i:S^1\to S^1$ given by
\begin{equation}
g_0(z)=z,\;\; g_1(z)=-z,\;\; g_2(z)=\bar z,\;\; g_3(z)=-\bar z,
\end{equation}
for all $z\in S^1$. We can identify $S^1$ with $\mathbb R/ \mathbb Z$ and for each ${[t]\in \mathbb R/ \mathbb Z}$  we have
\begin{equation}\label{wzor23}
g_0([t])=[t],\;\; g_1([t])=\left[t+\frac12\right],\;\; g_2([t])=[-t],\;\; g_3([t])=\left[-t+\frac12\right].
\end{equation}
Let ${\it{D}}=\langle g_i:i=0,1,2,3\rangle$. Then ${\it{D}}\cong\mathbb Z_2\times\mathbb Z_2$ and $g_3=g_1g_2$. We define an action ${\it{D}}^n$ on $T^n$ by
\begin{equation}
(t_1,\ldots,t_n)(z_1,\ldots,z_n)=(t_1z_1,\ldots,t_nz_n)
\end{equation}
for $(t_1,\ldots,t_n)\in{\it{D}}^n$ and $(z_1,\ldots,z_n)\in T^n=\underbrace{S^1\times\ldots\times S^1}_{n}$.

\noindent Any subgroup $\mathbb Z_2^d\sqsubseteq {\it{D}}^n$ defines $(d\times n)-$matrix with entries in ${\it{D}}$ which defines a matrix $P$ with entries in the set $S=\{0,1,2,3\}$
under the identification $i\leftrightarrow g_i$, $i=0,1,2,3$.

We have the following characterisation of the action of $\mathbb Z_2^d$ on $T^n$ and the associated orbit space $T^n/\mathbb Z_2^d$ via the matrix $P$.
Let $\mathbb Z_2^d\subseteq {\it{D}}^n$ and $P\in S^{d\times n}$. Then the action of $\mathbb Z_2^d$ on $T^n$ is free if and only if there is 1 in the sum of any distinct collection of rows of $P$. Group $\mathbb Z_2^d$ is the holonomy group of $T^n/\mathbb Z_2^d$ if and only if there is either 2 or 3 in each row of $P$.

Let us consider the epimorphisms $\alpha,\beta:{\it{D}}\to\mathbb F_2=\{0,1\}$ where values of $\alpha$ and $\beta$ on ${\it{D}}$ are given by

\begin{center}
	\begin{tabular}{|c|c|c|c|c|} \hline
		&0&1&2&3 \\ \hline\hline
		$\alpha$&0&1&1&0\\ \hline
		$\beta$&0&1&0&1\\ \hline 
	\end{tabular}\;.
\end{center}
For $j=1,2,\ldots,n$ and $\mathbb Z_2^d\subseteq{\it{D}}^n$ we define epimorphisms
\begin{equation}
\alpha_j:\mathbb Z_2^d\subseteq{\it{D}}^n\xrightarrow{pr_j} {\it{D}}\xrightarrow{\alpha}\mathbb F_2,\;\;
\beta_j:\mathbb Z_2^d\subseteq{\it{D}}^n\xrightarrow{pr_j} {\it{D}}\xrightarrow{\beta}\mathbb F_2
\end{equation}
by
$
\alpha_j(t_1,\ldots,t_n)=\alpha(t_j),\;\;
\beta_j(t_1,\ldots,t_n)=\beta(t_j).
$
Using definitions of $\alpha$ and $\beta$ and the translations given by (\ref{wzor23}), we obtain the following lemma.
\begin{lm}{\rm(\cite{LPPS})}
	Suppose a subgroup $\mathbb Z_2^d\subseteq{\it{D}}$ acts free on $T^n$. Then a holonomy representation $\varphi:\mathbb Z_2^d\to \operatorname{GL}(n,\mathbb Z)$ of the flat manifold $T^n/\mathbb Z_2^d$
	is given by
	$$\varphi(x)=\operatorname{diag}\left[(-1)^{(\alpha_1+\beta_1)(x)},\ldots,(-1)^{(\alpha_n+\beta_n)(x)}\right],$$
	for all $x\in\mathbb Z_2^d$.
\end{lm}
Since $H^1(\mathbb Z_2^d,\mathbb F_2)=\operatorname{Hom}(\mathbb Z_2^d,\mathbb Z_2)$ we can view $\alpha_i$ and $\beta_i$ as 1-cocycles and define
\begin{equation}
\label{wzor33}
\theta_j=\alpha_j\cup \beta_j\in H^2(\mathbb Z_2^d,\mathbb F_2),\end{equation}
where $\cup$ denotes the cup product. It is well known that
$H^*(\mathbb Z_2^d,\mathbb F_2)\cong \mathbb F_2[x_1,\ldots,x_d]$ where $\{x_1,\ldots,x_d\}$ is a basis of $H^1(\mathbb Z_2^d,\mathbb F_2)$. Hence, elements $\alpha_j$ and $\beta_j$ correspond to
\begin{equation}\label{wzor34}
\alpha_j=\sum_{i=1}^d\alpha(pr_j(b_i))x_i,\;\;\; \beta_j=\sum_{i=1}^d \beta(pr_j(b_i))x_i\in \mathbb Z_2^d[x_1,\ldots,x_d],
\end{equation}
where $\{x_1,\ldots,x_d\}$ is the standard basis of $\mathbb Z_2^d$ and $j=1,\ldots,n$ (\cite{CMR}, Proposition 1.3). Moreover, from the definition of matrix $P\in S^{d\times n}$
we can write equations (\ref{wzor33}) and (\ref{wzor34}) as follows
\begin{equation}\label{wzor35}
\alpha_j=\sum_{i=1}^d \alpha(P_{ij})x_i,\;\; \beta_j=\sum_{i=1}^d \beta(P_{ij})x_i,\;\;\theta_j^P=\alpha_j\beta_j.
\end{equation}
There is an exact sequence
$$0\to H^1(\mathbb Z_2^d,\mathbb F_2)\xrightarrow{\pi^*}H^1(\Gamma,\mathbb F_2)\xrightarrow{\iota^*}H^1(\mathbb Z^n,\mathbb F_2)\xrightarrow{d_2}
H^2(\mathbb Z_2^d,\mathbb F_2)\xrightarrow{\pi^*}H^2(\Gamma,\mathbb F_2)$$
where $d_2$ is the transgression and $\pi^*$ is induced by the quotient map $\pi:\Gamma\to \mathbb Z_2^d$, \cite{E}.
\begin{prop}{\rm(\cite{LPPS})}\label{Prop32}
	Suppose $\mathbb Z_2^d$ acts freely and diagonally on $T^n$. Let $M=T^n/\mathbb Z_2^n$, $\Gamma=\pi_1(M)$ and consider the associated to the group extension of (\ref{wzor21}). Then
	\begin{enumerate}
		\item $\forall_{1\leq l\leq n\;\;}\theta_l=d_2(\varepsilon_l)$, where $\{\varepsilon_1,\ldots,\varepsilon_n\}$is the basis of $H^1(\mathbb Z^n,\mathbb F_2)$ dual to the standard basis of $\mathbb Z^n\otimes\mathbb F_2$,
		\item the total Stiefel-Whitney class of $M$ is
		$$w(M)=\pi^*\left(\prod_{j=1}^{n}(1+\alpha_j+\beta_j)\right)\in H^*(\Gamma,\mathbb F_2)=H^*(M,\mathbb F_2).$$
	\end{enumerate}
\end{prop}

One can see by part (1) of Proposition \ref{Prop32} that the image of differential $d_2$ is an ideal generated by $\theta_j$ and
$$\langle\operatorname{Im}(d_2)\rangle=\langle \theta_1,\ldots,\theta_n\rangle\subseteq\mathbb F_2[x_1,x_2,\ldots,x_n.]$$
For matrix $P\in S^{d\times n}$, using (\ref{wzor35}), we set $\operatorname{Id}_P=\langle \theta_1^P,\ldots,\theta_n^P\rangle$ and we call this the {\it characteristic ideal} of $P$. The quotient $C_P=\mathbb F_2[x_1,\ldots,x_d]/\operatorname{{Id}}_P$ we call {\it characteristic algebra} of $P$.
\begin{cor}{\rm(\cite{LPPS})}
	Suppose $\mathbb Z_2^d$ acts freely and diagonally on $T^n$. There is a canonical homomorphism of graded algebras $\Phi:C\to H^*(T^n,\mathbb Z_2^d,\mathbb F_2)$ such that $\Phi([w])=w(T^n/\mathbb Z_2^d)$ where $[w]$ is the class of
	\begin{equation}\label{wzor37}
	w=\prod_{j=1}^n\left(1+\alpha_j+\beta_j\right)\in\mathbb F_2[x_1,\ldots,x_d].
	\end{equation}
	Moreover, $\Phi$ is a monomorphism in degree less that or equal to two.
\end{cor}
\begin{defi}
	Given a matrix $P\in S^{d\times n}$, we define the Stiefel-Whitney class of $P$, to be the class $[w]\in C_P$ defined by (\ref{wzor37}).
\end{defi}
\begin{cor}{\rm (\cite{LPPS})}
	Suppose $P\in S^{d\times n}$ is free and $T^n/\mathbb Z_2^d$ is the corresponding flat manifold. Then $\Phi(w(P))=w(T^n/\mathbb Z_2^d)$.
\end{cor}

Now, we describe a real Bott manifold $M(A)$. Let $A=[a_{ij}]$ be an strictly upper triangular matrix with entries 0 or 1 and let $s_i$, $i=1,\ldots,n$ be Euclidean motions on $\mathbb R^n$ defined by

\begin{equation}\label{gener}
s_i=\left(\operatorname{diag}\left[1,\ldots,1,(-1)^{a_{i,i+1}},\ldots,(-1)^{a_{i,n}}\right],\left(0,\ldots,0,\frac12,0\ldots,0\right)^T\right)
\end{equation}
where $(-1)^{a_{i,i+1}}$ is at the $(i+1, i+1)$ position and
$\frac{1}{2}$ is the $i-$th coordinate of the column, $i =
1,2,...,n-1.$ $s_{n} =
\left(I,\left(0,0,...,0,\frac12\right)\right)\in E(n).$ The group $\Gamma(A)$ generated by $s_1,\ldots,s_n$ is crystallographic group. The subgroup generated by
$s_{1}^{2},s_{2}^{2},...,s_{n}^{2}$ consists of all transitions by $\mathbb Z^n$. The action of $\Gamma(A)$ on $\mathbb R^n$ is free and the orbit space $\mathbb R^n/\Gamma(A)$ is compact.

\section{Main results}

We know that $s_i$ (\cite{KM}) are generators of the crystallographic group $\Gamma(A)$. Using the same methods as in \cite{PS}, \cite{LPPS}, for each strictly upper triangular matrix $A=[a_{ij}]$ which generates the fundamental group of real Bott manifold $M(A)$ we have $P-$matrix $P_A=[p_{ij}]$ with diagonal entries 1,  entries 0 or 2 in the upper triangular part and entries 0 in the lower triangular part.
Since $\beta(0)=\beta(2)=1$ and $\beta(1)=1$ we get
$$\beta_j=\sum_{i=1}^d\beta\left(p_{ij}\right)x_i=\beta\left(p_{jj}\right)x_j=x_j \text{   for all }1\leq i \leq n.$$
Now, $\alpha(0)=0$, $\alpha(1)=\alpha(2)=1$ and
$$\begin{aligned}
\alpha_j+\beta_j&=\sum_{i=1}^{j-1}\alpha\left(p_{ij}\right)x_i+\sum_{i=1}^{j-1}\beta(p_{ij})x_i
=\sum_{\substack{i=1, i\ne j}}^{j-1}\alpha(p_{ij})x_i+\sum_{\substack{i=1, i=j}}^{j-1}\left(\alpha(p_{ij})+\beta(p_{ij})\right)x_i\\
&=\sum_{\substack{i=1, i\ne j}}^{j-1}\alpha(p_{ij})x_i+\left(\alpha(p_{jj})+\beta(p_{jj})\right)x_j
=\sum_{i\ne j}\alpha(p_{ij})x_i=\sum_{i=1}^{j-1}b_{ji}x_i,
\end{aligned}$$
where $b_{ji}\in\{0,1\}$ for all $i=1,2,\ldots,j-1$. So, we get
\begin{equation}\alpha_i+\beta_i=\sum_{j=1}^{i-1}b_{ij}x_j,\label{wzor2}\end{equation}
where $b_{ij}\in\{0,1\}$ for all $j=1,2,\ldots,i-1$.

Now, let us consider a real Bott manifold with a K\"{a}hler structure. We will denote this manifold by RBK manifold. Let the column $P^k$  of matrix $P_A$ has $p_{kk}=1$ and all others entries equal to 0. Then
$$\alpha_k=x_k=\beta_k,$$ so $x_k^2\in \operatorname{Id}_{P_A}$.

\begin{theo}\label{theo1}
	Let $P_A$ be a $P-$matrix of $2n-$dimensional RBK manifold $M(A)$. Then
	$$w_2(M(A))=\sum_{i=1}^{n}\left(b_{j_1i}+b_{j_2i}+\ldots+ b_{j_ki}\right)x_i^2\in \mathbb F_2[x_1,x_2,\ldots,x_n]$$
	and  $M(A)$ has a Spin-structure if and only if either $\left(b_{j_1i}+b_{j_2i}+\ldots+ b_{j_ki}\right)=0\mod 2$ or $x_i^2\in\operatorname{Id}_{P_A}$ for all $i=1,\ldots,n$.
\end{theo}

\noindent{\rm Proof.} $M(A)$ is $2n-$dimensional RBK manifold, so from (\ref{wzor2})
we get
$$\alpha_{j_{k+n}}+\beta_{j_{k+n}}=\alpha_{j_k}+\beta_{j_k}=\sum_{i=1}^{n}b_{j_{k}i}x_i,$$
where $b_{j_ki}\in\{0,1\}$ for all $i=1,2,\ldots,j_k-1$ and
$$\begin{aligned}
w(M(A))&=\prod_{k=1}^{2n}\left(1+\alpha_{j_k}+\beta_{j_k}\right)=\prod_{k=1}^{n}\left(1+\alpha_{j_k}+\beta_{j_k}\right)^2\\
&=\prod_{k=1}^{n}\left(1+\sum_{i=1}^{j_k-1}b_{j_ki}x_i\right)^2=\prod_{k=1}^{n}\left(1+\sum_{i=1}^{j_k-1}b_{j_ki}x_i^2\right)\\
&=1+\sum_{i=1}^{n}\left(b_{j_1i}+b_{j_2i}+\ldots+ b_{j_ni}\right)x_i^2+\ldots \;.
\end{aligned}$$
From the above consideration we have
$$w_2(M(A))=\sum_{i=1}^{n}\left(b_{j_1i}+b_{j_2i}+\ldots+ b_{j_ni}\right)x_i^2.$$
It is well known that the manifold $M(A)$ has the Spin-structure if $w_2(M(A))=0$. In our case
$$w_2(M(A))=\sum_{i=1}^{n}\left(b_{j_1i}+b_{j_2i}+\ldots+ b_{j_ni}\right)x_i^2=0$$
if either $\left(b_{j_1i}+b_{j_2i}+\ldots+ b_{j_ni}\right)=0\mod 2$ or $x_i^2\in\operatorname{Id}$ for all $i=1,\ldots,n$.

\hskip13cm$\square$

At the end we consider a special case of RBK manifold.

\begin{lm}\label{lemma2}
	Let $M(A)$ be a RBK manifold with matrix $P$. Let $k$ be an even number and let $P^{i_1}=P^{i_2}=\ldots=P^{i_{2k}}$ be columns with nonzero entries and all others columns of matrix $P$ have only 0 entries. Then the manifold $M(A)$ has the Spin-structure.
\end{lm}
\noindent
{\rm Proof.} Since $P^{i_1}=P^{i_2}=\ldots=P^{i_{2k}}$ we get
$$\alpha_{i_{1}}+\beta_{i_{1}}=\ldots=\alpha_{i_{2k}}+\beta_{i_{2k}}=\sum_{k=1}^{i_1}b_{i_1k}x_k,$$
where $b_{i_1k}\in\{0,1\}$ for all $k=1,2,\ldots,i_1$ and
from the proof of Theorem \ref{theo1} and since $k=2s$ we get
$$\begin{aligned}
w(M(A))&=\prod_{l=1}^{2k}\left(1+\alpha_{j_l}+\beta_{j_l}\right)=\left(1+\alpha_{j_l}+\beta_{j_l}\right)^{2k}=\left(1+\sum_{r=1}^{j_l-1}b_{j_lr}x_r\right)^{2k}\\
&=\left(1+\sum_{r=1}^{j_l-1}b_{j_lr}x_r^2\right)^{k}=\left(1+\sum_{r=1}^{j_l-1}b_{j_lr}x_r^2\right)^{2s}=\left(1+\sum_{r=1}^{j_l-1}b_{j_lr}x_r^4\right)^{s}.
\end{aligned}$$
So,
$w_2(M(A))=0$
and RBK manifold has the Spin-structure.

\hskip13cm$\square$

\begin{ex}
	Let $$A=\left[\begin{matrix}1&0&2&2&2&2\\0&1&2&2&2&2\\0&0&1&0&2&2\\0&0&0&1&2&2\\0&0&0&0&1&0\\0&0&0&0&0&1\end{matrix}\right]$$
	be a matrix of the manifold $M(A)$. Then
	$$\begin{aligned}\operatorname{Id}_A&=\left\{x_1^2,x_2^2,x_1x_3+x_2x_3+x_3^2,x_1x_4+x_2x_4+x_4^2,\right.\\ &\left.x_1x_5+x_2x_5+x_2x_6+x_4x_5+x_5^2,
	x_1x_6+x_2x_6+x_3x_6+x_4x_6+x_6^2\right\}\\
	\end{aligned}$$
	and
	$$\begin{aligned}w(M(A))=x_3^2+x_4^2.
	\end{aligned}$$
	Since $x_3^2,x_4^2\not\in\operatorname{ Id}_A$, so
	$w_2(M(A))=x_3^2+x_4^2\ne 0$
	and $M(A)$ has no Spin-structure.
\end{ex}

Let $A$ be the matrix of $2n-$dimensional  RBK manifold $M(A)$ and let
$S_k=\sum_{i=1}^na_{kj_i}$ for all $1\leq k\leq 2n$. Then from Theorem \ref{thI} and following Theorem \ref{theo1} we get

\begin{cor}
	Let $A$ be a matrix of a $2n-$dimensional RBK manifold $M(A)$. The manifold $M(A)$ has the Spin-structure if or $S_k=0\mod 2$ or  $S_k=1\mod 2$ and the column $A^k$ has only entries 0, for all $1\leq k\leq 2n$.
\end{cor}
\begin{ex}
	Let $$A=\left[\begin{matrix}0&0&1&1&1&1\\0&0&1&1&1&1\\0&0&0&0&1&1\\0&0&0&0&1&1\\0&0&0&0&0&0\\0&0&0&0&0&0\end{matrix}\right]$$
	be a matrix of a manifold $M(A)$. Then $S_1=S_2=2\mod 2=0, S_3=S_4=1\mod 2=1,	S_5=S_6=0\mod 2=0$ and there are entries equal to 1 in columns $A^3$ and $A^4$, so $M(A)$ has no Spin-structure.
\end{ex}

\vskip 2mm
\noindent
Maria Curie-Sk{\l}odowska University,\\
Institute of Mathematics\\
pl. Marii Curie-Sk{\l}odowskiej 1\\
20-031 Lublin, Poland\\
E-mail: anna.gasior@poczta.umcs.lublin.pl

\end{document}